\documentclass[12pt]{amsart}
\usepackage{amssymb}
\usepackage{amsthm}
\usepackage{mathrsfs}
\usepackage{epsfig}

\newtheorem{thm}{Theorem}[section]

\newtheorem{cor}[thm]{Corollary}
\newcommand{\beg}[2]{\begin{equation}\label{#1}#2\end{equation}}

\newcommand{\homotopic}{\simeq}							

\newcommand{\rref}[1]{(\ref{#1})}
\newcommand{\cA}{\mathcal{A}}

\newcommand{\cF}{\mathcal{F}}

\newcommand{\Z}{\mathbb{Z}}

\title[$(\Z/p)^n$-Equivariant Fixed Points of $H\Z/p$]{On the Coefficients of $(\Z/p)^n$-Equivariant Ordinary
Cohomology with Coefficients in $\Z/p$}
\author{John Holler and Igor Kriz}

\begin{document}
\maketitle

\begin{abstract}
This note contains a generalization to $p>2$ of the authors' previous calculations \cite{hk} of the coefficients
of $(\Z/2)^n$-equivari-ant ordinary cohomology with coefficients in the constant $\Z/2$-Mackey
functor. The algberaic results by S.Kriz \cite{sk} allow us to calculate the coefficients of the geometric
fixed point spectrum $\Phi^{(\Z/p)^n}H\Z/p$, and more generally, the $\Z$-graded coefficients
of the localization of $H\Z/p_{(\Z/p)^n}$ by inverting any chosen set of embeddings $S^0\rightarrow S^{\alpha_i}$
where $\alpha_i$ are non-trivial irreducible representations. We also calculate the $RO((\Z/p)^n)^+$-graded 
coefficients of $H\Z/p$, which means the cohomology of a point indexed by an actual (not 
virtual) representation. (This is the ``non-derived" part, which has a nice algebraic description.)
\end{abstract}

\section{Introduction}

In \cite{hk}, we
made progress on calculating the $RO(G)$-graded fixed points of $H\Z/2_G$ where
$G=(\Z/2)^n$. Here $H\Z/2_G$ denotes the ``ordinary" $RO(G)$-(co)homology theory
corresponding to the ``constant" Mackey functor $\Z/2$ (see \cite{lmm}).
More precisely, we calculated an explicit answer in (homologically graded) dimensions $*-V$ where 
$*\in \Z$, and $V$ is a finite-dimensional real $G$-representation with $V^G=0$. In other dimensions,
we gave a ``reasonably small" chain complex calculating the answer, but the answer is chaotic, 
and difficult to present explicitly.

The most interesting 
part of the calculation \cite{hk} was the calculation of the ``geometric fixed points" of $H\Z/2_G$ (\cite{lms},
Section II.9). Here, an algebraic structure appeared which we were able to calculate by hand, but did not understand,
and did not know how to generalize to an odd prime.

The purpose of this note is to generalize the calculations of \cite{hk} to an odd prime $p$, with a particular emphasis
on the geometric fixed points of $H\Z/p_G$ where $G=(\Z/p)^n$. The crucial development since the paper 
\cite{hk} was written was the appearance of the paper \cite{sk} by Sophie Kriz, who identified the algebraic structure
we found as the {\em reciprocal space} of a hyperplane arrangement in the case of all non-trivial hyperplanes
in the $n$-dimensional vector space over $\Z/2$, and also investigated the appropriate algebraic analogue for
an odd prime, which she called the {\em super-reciprocal space}. 

With the algebraic observations of \cite{sk} in hand, we were able to completely understand the geometric
fixed points of $H\Z/p_G$ in the case of an odd prime $p$ (and also, the treatment of \cite{hk} greatly simplifies).
The results of \cite{sk} also suggest that the (super)reciprocal spaces of other hyperplane arrangments should
appear as $\Z$-graded coefficients of $G$-equivariant spetra, and with the present techniques, this is, in fact, 
easy to confirm. The purpose of this note is to record these developments.

\section{Geometric Fixed Points}

A ''cube" spectral sequence was developed in \cite{hk}, and was used to give a complete characterization of the ring structure of the geometric fixed points of $H\Z/p$ for $p=2$. In the present note, we calculate the coefficients of the geometric fixed points of the ordinary equivariant cohomology of the constant $\Z/p$ Mackey functor for $G = (\Z/p)^n$. To this end we will make use of the same spectral sequence.

Let $G = (\Z/p)^n$, where $p >2$. Following the notation from [2], $\cF[H]$ is the family of all subgroups $K \subset G$ such that $H$ is not contained in $K$. $\cF(H)$ is the family of all subgroups $K \subset H$. Given any family of subgroups $\cF$, there is an ''$\cF$-universal space", which is a $G$-CW complex $E\cF$ satisfying

$$
 E\cF^K =
  \begin{cases} 
      *    & \text{  if  $K \in \cF$} \\
     \emptyset & \text{ if $K \notin \cF$} \\
  \end{cases}
$$
$E\cF(H)$ is often written as $EG/H$ since as a $G/H$-space it is the universal space associated to the family containing only the trivial subgroup. For any family of subgroups of $G$, we have the isotropy separation sequence
$$
E\cF_+ \rightarrow S^0 \rightarrow \widetilde{E\cF}
$$
For a G-spectrum $X$, the geometric fixed points are given by 
$$\Phi^G(X) = (\widetilde{E\cF[G]} \wedge X)^G,$$ 
which is the homotopy cofiber of $E\cF[G]_+ \wedge X \rightarrow X$. There is a particularly nice model for $\widetilde{E\cF[G]}$ which allows us to view the geometric fixed points of a spectrum as an iterated cofiber - namely $\widetilde{E\cF[G]} = S^{\oplus_{V}\infty V}$ where the direct sum is taken over all irreducible (real) representations of $G$, and the $\infty$ means we are taking the diret limit over finitely many copies of each $V$. Given a maximal subgroup $H \subset G$, $ S^{\oplus_{V}\infty V}$ is a model for $\widetilde{EG/H}$ where now the direct sum is taken over irreducible representations of $G$ that are $H$-trivial. Each of these representations corresponds to an irreducible $G/H = \Z/p$ representation, and so there are exactly $p-1$ of them, each of dimension 2. \newline
\indent Using these models, it is clear that the $\Phi^G(H\Z/p)$ is given as the colimit of the diagram formed by smashing together the maps
$$
EG/H_+ \rightarrow S^0
$$
over all maximal subgroups $H \in G$, smashing with $H\Z/p$, and then taking $G$-fixed points. Each maximal subgroup $H$ can be realized uniquely up to scalar multiplication as the kernel of a map $\alpha: G \rightarrow \Z/p$. For each $H$ we choose the unique corresponding $\alpha$ of the form $(*, \ldots, *,1,0,\ldots,0)$ where the $*$'s can be any elements of $\Z/p$. \newline
\indent By examining fixed points it's easy to see that for any collection of subgroups $H_1, \ldots, H_k$ of $G$,
$$
EG/H_1 \times \cdots \times EG/H_k \homotopic EG/(H_1 \cap \cdots \cap H_k)
$$
Therefore, taking the first quadrant spectral sequence associated to this cube diagram, we get the following $E^1$-page:
\beg{espec1}{\resizebox{\textwidth}{!}{%
 $E^1_{*,s}=
  \begin{cases} 
     ( H\Z/p)_*   & \text{  if  s = 0} \\
     \bigoplus_{S \in A_s} Sym(G/(\cap_{\alpha \in S}ker(\alpha))^*) \otimes \bigwedge (G/(\cap_{\alpha \in S}ker(\alpha))^*)\cdot y_S& \text{ if s $>$ 0} \\
  \end{cases}$%
}
}
where $A_s$ is the collection of cardinality $s$ subsets of the set $\cA=G^*\smallsetminus \{0\}$. 
Unlike the case $p=2$, the polynomial generators here are in dimension 2, and the exterior generators (which do not show up at all in the $p=2$ case), are in dimension 1.

\bigskip
With this particular choice of $\alpha$'s, the description of the $E^2$-page given in \cite{hk} can be imported here in its analogous form. Explicitly,
$$
E^2 = \bigoplus_{y_S \in F_n} Sym((G/ \cap(ker(\alpha) | \alpha \in S))^*) \otimes \bigwedge((G/ \cap(ker(\alpha) | \alpha \in S))^*) \cdot y_S
$$
where $F_n$ is given in Lemma 1 of \cite{hk} (with 2 replaced by $p$). Explicitly,
$$F_1=\{y_\emptyset, y_{\{(1)\}}\},$$
$$F_n=F_{n-1}\cup \{y_{S\cup \{x\}}\mid S\in F_{n-1}, x\in (\Z/p)^{n-1}\times\{1\}\}.$$
In other words, $F_n$ consists of the basis elements $y_S$ where $S$ are all $\Z/p$-linearly independent (in $G^*$)
subsets in (not necessarily reduced) row echelon form with respect to reversed order of columns.

We maintain the notation in \cite{hk} by denoting by $y_S$ the generator corresponding to $S$. 

\vspace{3mm}
We next consider the $G$-spectrum
\beg{elangle}{\langle H\Z/p_G\rangle=\widetilde{E}\cF[G] \wedge F(EG_+, H\Z/p).}
We have, of course, a canonical map 
\beg{elangle1}{\widetilde{E}\cF[G]\wedge H\Z/p_G\rightarrow \langle H\Z/p_G\rangle.
}
For the coefficients of the Borel cohomology spectrum, we may write
(following the notation of \cite{sk})
\beg{eborelc}{F(EG_+, H\Z/p)_*=Z/p[x_{\alpha_1},\dots,x_{\alpha_n}]\otimes \Lambda_{\Z/p}[dx_{\alpha_1},
\dots,dx_{\alpha_n}]}
where $\alpha_1,\dots,\alpha_n$ is a $\Z/p$-basis of $G^*$, and the generators $x_\alpha$, $dx_\alpha$ have 
homological dimension $-2, -1$, respectively.

The Borel cohomology spectrum $F(EG_+, H\Z/p)$ is periodic with respect to smashing with 
any virtual sphere $S^\alpha$ where $\alpha$ is in the augmentation ideal of $RO(G)$, so the coefficients
of $\langle H\Z/p_G\rangle$ can be computed from the coefficients of Borel cohomology 
simply by inverting Euler classes, which are all the $p^n-1$
non-zero $\Z/p$-linear combinations $z_\alpha$ of $x_{\alpha_i}$:
\beg{eborelc1}{\langle H\Z/p_G\rangle_*=(\prod_\alpha z_\alpha)^{-1}(
Z/p[x_{\alpha_1},\dots,x_{\alpha_n}]\otimes \Lambda_{\Z/p}[dx_{\alpha_1},
\dots,dx_{\alpha_n}]).
}
Following \cite{sk}, one denotes $t_\alpha=(z_\alpha)^{-1}$.
However, we can also compute \rref{eborelc1} by a spectral sequence analogous to \rref{espec1}. The 
$E^1_{*,s}$-term of that spectral sequence is \rref{eborelc} for $s=0$, and for $s>0$, using a notation
analogous to \rref{espec1}, the $\Z/p$-module
attached to the generator $y_S$ is
\beg{etate1}{
(Sym(G/(\cap_{\alpha \in S}ker(\alpha))^*) \otimes \bigwedge (G/(\cap_{\alpha \in S}ker(\alpha))^*))\otimes 
\Z/p[x_{\alpha_i}]\otimes \Lambda[dx_{\alpha_i}]
}
where $\alpha_i$ is the set of free generators of the $\Z/p$-vector space 
$$G^*/\langle S\rangle=(\bigcup_{\alpha\in S} Ker(\alpha))^*.$$
The $E_2$-term,
again, is the sum of the terms with $y_S\in F_n$. Now the generator $y_S$ is a permanent cycle, representing
$$\prod_{\alpha\in S} t_\alpha,$$
while the polynomial generators in \rref{etate1} are $t_\alpha$, $\alpha\in S$, and the exterior generators
are $u_\alpha=t_\alpha dz_\alpha$, $\alpha\in S$.
Thus, this spectral sequence collapses to $E^2$. Since the spectral sequence \rref{espec1} naturally maps to it
by the map \rref{elangle1}, and the map of spectral sequences is an inclusion on $E^1$ and $E^2$, we see that the spectral sequence \rref{espec1} also collapses to $E^2$ and that moreover, $\Phi^G H\Z/p$ is by the inclusion
induced by \rref{elangle1} on coefficients identified with the subring of \rref{eborelc1} generated by $t_\alpha$, 
$u_\alpha$ for $\alpha \in G^*\smallsetminus \{0\}$. Using Theorem 6 of \cite{sk} to all 
``hyperplanes through the origin" in $(\Z/p)^n$, we then obtain the following

\begin{thm}\label{egeom}
The ring $\Phi^G(H\Z/p)_*$ is isomorphic to the quotient of 
$$\Z/p[t_\alpha\mid \alpha\in G^*\smallsetminus\{0\}]\otimes 
\Lambda_{\Z/p}[u_\alpha\mid \alpha\in G^*\smallsetminus\{0\}]$$
modulo the ideal generated by
$$t_\alpha-kt_{k\alpha},\;k\in \Z/p^\times,$$
$$t_\alpha t_\beta+t_\beta t_\gamma+t_\gamma t_\alpha,$$
$$u_\alpha t_\beta+u_\alpha t_\gamma-u_\beta t_\gamma -u_\gamma t_\beta,$$
$$u_\alpha u_\beta+u_\beta u_\gamma+u_\gamma u_\alpha,\; \alpha+\beta+\gamma=0,$$
where $\alpha, \beta,\gamma\in G^*\smallsetminus \{0\},\;\; \alpha+\beta+\gamma=0$.
\end{thm}
\qed

\begin{cor} The Poincare series of $(\Phi^GH\Z/p)_*$ is
$$
\frac{1}{(1-x)^n}\prod_{i=1}^{n}(1+(p^{i-1}-1)x)
$$
\end{cor}
\qed

\section{$RO(G)^+$-graded coefficients and other localizations}

The method of Section 4 of \cite{hk} now also applies verbatim, calculating the ``negative part" of the
$RO(G)$-graded coefficient ring of $H\Z/p_G$. Let, for a non-trviial $2$-dimensional representation 
$\alpha$ of $G$, $a_\alpha\in \pi_{-\alpha}^G(S)$ be the class represented by the inclusion
$$S^0\subset S^\alpha.$$

\begin{thm}
\label{thm2}
The ring $R_{p,n}=(H\Z/p_G)_{*-V}$ where $*$ 
stands for an integer, and $V$ for a $G$-representation with $V^G=0$ is the
subring of 
$$\Phi^G(H\Z/p)_*[a_\alpha \mid \alpha\in G^*\smallsetminus \{0\}]$$
generated by 
$$a_\alpha t_\alpha,\; a_\alpha u_\alpha.
$$
\end{thm} 
\qed

The Poincar\'{e} series calculation of Theorem 5 of \cite{hk} also translates 
verbatim from \cite{hk} with $2$ replaced by $p$ and $m_\alpha$ replaced by $2m_\alpha$.

By choosing a set $S\subseteq G^*\smallsetminus \{0\}$, the $\Z$-graded part $Q_{p,n,S}$ of the ring
$$(\prod_{\alpha\in S} a_\alpha)^{-1} R_{p,n}$$
is the $\Z$-graded coefficient ring of the $G$-spectrum
$$\bigwedge_{\alpha\in S} S^{\infty \alpha}\wedge H\Z/p_G.$$
From Theorem \ref{thm2}, it follows that this is the subring of \rref{eborelc1} generated by 
$t_\alpha, u_\alpha$, $\alpha\in S$. Since localization by a non-zero divisor is an injective map, this is the 
same as the subring of 
$$
(\prod_{\alpha\in S} z_\alpha)^{-1}(
Z/p[x_{\alpha_1},\dots,x_{\alpha_n}]\otimes \Lambda_{\Z/p}[dx_{\alpha_1},
\dots,dx_{\alpha_n}])
$$
generated by $t_\alpha, u_\alpha$. These rings were also calculated in \cite{sk}, Theorem 6.

\end{document}